\documentclass[12pt]{article}

\usepackage[bottom]{footmisc}

\usepackage{amsmath}
\usepackage{latexsym}
\usepackage{amssymb}

\newcommand{\CC}{{\cal C}}

\newtheorem{thm}{Theorem}

\newtheorem{Defn}[thm]{Definition}
\newtheorem{Remark}[thm]{Remark}
\newtheorem{Note}[thm]{Note}

\newtheorem{Example}[thm]{Example}
\newtheorem{Examples}[thm]{Examples}
\newtheorem{Problems}[thm]{Problems}

\newtheorem{Problem}[thm]{Problem}
\newtheorem{Notation}[thm]{Notation}
\newtheorem{Number}[thm]{\!\!}

                  {\nopagebreak\hspace*{\fill}$\Box$\medskip\medskip\par}

\newcommand{\n}{\rm}

\newcommand{\R}{{\mathbb R}}

\newcommand{\Z}{{\mathbb Z}}
\newcommand{\C}{{\mathbb C}}
\newcommand{\K}{{\mathbb K}}

\newcommand{\g}{{\mathfrak g}}

\newcommand{\Aut}{\mbox{\n Aut}}

\newcommand{\pr}{\mbox{\rm pr}}

\newcommand{\id}{\mbox{\n id}}

\newcommand{\eps}{{\varepsilon}}

\newcommand{\msk}{\medskip}
\newcommand{\ssk}{\smallskip}
\newcommand{\nin}{\noindent}

\newcommand{\ttt}{{\mathbf t}}
\newcommand{\xxx}{{\mathbf x}}

\begin{document}

\title{Difference Problems and Differential Problems} 

\author{Wolfgang  Bertram}

\maketitle

\footnotesize{
{\noindent{\bf Abstract.\/}
We state some elementary problems concerning the relation between
{\em difference calculus} and {\em differential calculus},
and we try to convince the reader that,
in spite of the simplicity of the statements, a solution of these
problems would be a significant contribution to the understanding
of the foundations of differential and integral calculus.
}

\bigskip
\noindent
{\bf AMS subject classification:}
14A25, 
39A12, 
39A13, 
46S10, 
58C20 

\bigskip
\noindent
{\bf Key words:} difference calculus, differential calculus, 
divided differences,
super differential calculus, Pansu calculus,
$q$-derivative, scalar extension

\vspace{10mm}
\nin
W. Bertram, 
Institut \'Elie Cartan Nancy, Nancy-Universit\'e, CNRS, INRIA, 
Boulevard des Aiguillettes, B.P. 239, F-54506 Vand\oe{}uvre-l\`{e}s-Nancy,
 France;
{\tt bertram@iecn.u-nancy.fr}

\section*{Introduction}

Differential calculus can be seen
as the ``continuous extension  of difference calculus to singular
parameters''.\footnote{I thank Wilhelm Kaup for pointing out to me that
this observation is folklore.} This simple observation
 has been used in \cite{BGN04} 
for developing a general approach to differential calculus in
arbitrary dimension and over very general base fields and rings:
let $\K$ be a commutative ring with unit $1$ (for a first reading, think 
of $\K=\R$ or $\K=\C$) und $V$, $W$ be $\K$-modules.
For any subset $U \subset V$ and any map $f:U \to W$ we define 
the {\em first order difference quotient} by
$$
f^{]1[} (x,v,t) := {f(x+tv) - f(x) \over t},
\eqno (0.1)
$$
wherever this makes sense, i.e., for $(x,v,t)$ belonging to the set
$$
U^{]1[} = \{ (x,v,t) \in U \times V \times \K^\times \vert \, x+tv \in U \}.
\eqno (0.2)
$$
The map $f^{]1[}:U^{]1[} \to W$ thus defined
is called the {\em first order difference quotient map}.
If we want to stress the dependence on the ring $\K$, we write also
$f^{]1[}_\K$ and $U^{]1[}_\K$.
In order to keep track of all three variables, it is also useful to
introduce the {\em extended difference quotient map}
$$
\Delta^{]1[} f: U^{]1[} \to W^{]1[} 
, \quad
(x,v,t) \mapsto \bigl( f(x), f^{]1[} (x,v,t) , t \bigr)
\eqno (0.3)
$$
which satisfies, as is immediately checked, the functorial relations
$$
\Delta^{]1[} (g \circ f) = \Delta^{]1[} g \circ \Delta^{]1[} f, \quad
\quad \Delta^{]1[} (\id_U) = \id_{U^{]1[}}.
\eqno (0.4)
$$

\msk
{\em Differential calculus} is the continuous extension to the singular
parameter value $t=0$. More precisely, let
$$
U^{[1]} = \{ (x,v,t) \in U \times V \times \K \vert \, x+tv \in U \} ;
\eqno (0.5)
$$
this set contains $U \times V \times \{ 0 \}$.
Assume that $\K$ is a topological ring such that the unit group
$\K^\times$ is open dense in
$\K$; assume moreover that $V$ and $W$ are topological $\K$-modules and that
$U$ is open $V$; then $U^{]1[}$ is dense in $U^{[1]}$, and both are open 
in $U \times V \times \K$.
We say that $f$  is {\em of class $\CC^1$} if $f^{]1[}$ admits
an extension to a {\em continuous} map $f^{[1]}:U^{[1]} \to W$ (equivalently,
if $\Delta^{]1[}$ admits an extension to a continuous map
$\Delta^{[1]} f: U^{[1]} \to W^{[1]}$).
By density of $\K^\times$ in $\K$, this extension is unique, if it exists,
and hence the {\em differential} and the {\em tangent map} are well-defined by
$$
df(x)v:=f^{[1]}(x,v,0), \quad \quad
Tf(x)v:=(f(x),df(x)v).
\eqno (0.6)
$$
Again by density, Relation (0.4) implies that $T(g \circ f)=Tg \circ Tf$,
i.e., the chain rule. In a similar way, all elementary rules of differential
calculus (for instance, the linearity of $df(x):V \to W$)
follow ``by density'' from simple difference calculus.
In \cite{BGN04}, it has been shown that, in all relevant cases,
these definitions are equivalent to more classical ones; in particular,
mappings between finite dimensional real vector spaces are $\CC^1$ in our
sense iff they are $\CC^1$ in the usual sense, and this can even be used
to simplify several arguments of usual multivariable analysis --
see \cite{Be08} for
an elementary account. On the other hand, this approach can be
further generalized beyond the context of differential calculus over
topological fields and rings  (cf.\ \cite{BGN04}).

\msk
All this may look fairly trivial, but that impression changes completely
if we now turn to {\em higher order difference and differential calculus}.
In fact, definitions are made such that they can easily be iterated; for
instance, $f$ is called {\em of class $\CC^2$} if it is $\CC^1$ and if
$f^{[1]}$ (or, equivalently, $\Delta^{[1]}f$) is again $\CC^1$,
giving rise to $f^{[2]}:=(f^{[1]})^{[1]}$ (resp.\
$\Delta^{[2]}f := \Delta^{[1]} (\Delta^{[1]}f)$), and so on.
Already the formula for $f^{]2[}$ is quite unpleasant (see the formula given
 below), and we are far from understanding  the
 general maps $f^{[k]}$ and $\Delta^{[k]}f$ (which have $2^{k+1} - 1$
arguments!). The reader will certainly agree that differential calculus
is much simpler than difference calculus because, in the limit,
 it only retains 
the ``invariant'' and most important features of the wild structure of
the latter. However, it remains important to study {\em how} this
limit is attained, especially if we have in mind the problem of
{\em integration} or {\em anti-derivative}, i.e., of going back from
differential data to local or even global ones. 
As we will see, the problem of understanding the precise way how
differential calculus in imbedded into difference calculus is highly
non-trivial.

\msk
In this note, I shall present some problems related to higher
order difference and differential calculus, most of them completely open
and, I hope, challenging specialists as well as non-specialists.
My impression is that some of these problems are of a  foundational
nature,  tieing together algebra, combinatorics, linear and affine
geometry with analysis, and that it would be very useful to better understand
what is going on here. 

\msk
%
I have presented some aspects of these problems at the Eighth International
Workshop on Differential Geometry and its Applications, held in Cluj-Napoca 
in August 2007, and I would like to thank the organizers for their
great hospitality.

\section{Difference problems}


\nin {\bf Problem 1: Explicit formulae.} 
{\em Find an explicit formula for the higher order difference quotient
map $f^{]k[}$, resp.\ for its extended version
$\Delta^{]k[} f$.}

\msk \nin {\bf Comments.}
For $k=2$,  we have the explicit formula
\begin{eqnarray*}
& & f^{]2[}((x_1,v_1,t_1),(x_2,v_2,t_2),t_3)   \cr
& & \quad \quad =
{f^{]1[}((x_1,v_1,t_1)+t_3 (x_2,v_2,t_2)) - f^{]1[}(x_1,v_1,t_1)
\over t_3} \cr
& & \quad \quad =
{f\bigl( x_1 + t_3 x_2 + (t_1+t_2 t_3) ( v_1 + t_3 v_2) \bigr) 
-f(x_1 + t_3 x_2)
\over  t_3 (t_1  + t_2 t_3) }
- {f(x_1+t_1 v_1) - f(x_1) \over t_1 t_3 } ,
\cr
\end{eqnarray*}
and the seven components of $\Delta^{]2[} f((x_1,v_1,t_1),(x_2,v_2,t_2),t_3)$
are
$$
\Big( \big(f(x_1),f^{]1[}(x_1,v_1,t_1),t_1\big),\big(f^{]1[}(x_1,x_2,t_3),
f^{]2[}((x_1,v_1,t_1),(x_2,v_2,t_2),t_3 \big), t_2),t_3 \Big).
$$
Clearly, it would be hopeless to try write out in this way
the formulae for general $k$. Thus we ask: 

\bigskip
\nin {\bf Problem 2: The structure.}
{\em Understand the structure of the higher order difference quotient
maps $f^{]k[}$, resp.\ of its extended version
$\Delta^{]k[} f$. Extend this structure (in the above mentioned
 setting of differential
calculus) to all singular parameter values.}

\msk \nin {\bf Comments.}
Both problems are related to each other, but are not equivalent.
The explicit formulae above don't really give an idea on
what goes on on the ``singular set'' (the set where 
some of the $t_i$'s are not invertible). Indeed, there exist 
non-trivial relations, such as
certain ``homogenity relations'' used in \cite{Gl04} (Appendix B) or
\cite{Gl05} (Lemma 6.8), or the following surprising relation from
 \cite{BGN04}, Chapter 5:
$$
f^{[2]}\big((x,v,0),(v,0,0),0 \big) =2f^{[2]}\big((x,v,0),(0,0,1),0 \big) \, .
$$
The explicit formula given above does not help at all to prove or
guess this relation. 
Yet, an explicit formula may turn out to be useful for attacking
Problem 2 -- for instance, in the above formula for $f^{]2[}$ we
see that for $t_2 = 0$ we get a rather familiar and more symmetric
expression.

\ssk
A first step towards the solution of
 Problems 1 and 2 consists in finding more appropriate
notation. The maps $f^{[k]}$ and $\Delta^{[k]}f$ have
 $2^{k+1}-1$ arguments: there are
$2^k$ ``space variables'' and $2^{k}-1$ ``scalar'' or ``(multi-)
time variables''.
For better book-keeping, we may write $u_0
 = (x_0,x_1,t_1)$ instead of $(x,v,t)$,
and use double indices on the next level; for instance
$$
(u_0,u_1,t_1)=
((x_{00},x_{01},t_{01}),(x_{10},x_{11},t_{11}),t_{10}),
$$
so that the formula for $f^{]2[}$ given above reads now
$$
 f^{]2[}\big((x_{00},x_{01},t_{01}),(x_{10},x_{11},t_{11}),t_{10}\big) =
$$
$$
{f\bigl( x_{00} + t_{10} x_{10} + (t_{01}+t_{11} t_{10}) 
(x_{01} + t_{10} x_{11}) \bigr) 
-f(x_{00} + t_{10} x_{10})
\over  t_{10} (t_{01}  + t_{11} t_{10}) }
- {f(x_{00}+t_{01} x_{01}) - f(x_{00}) \over t_{01} t_{10} } 
$$
and the seven components of $\Delta^{[2]}  f$ are 
(where we now put the multi-time variable $\ttt=(t_{01},t_{10},t_{11})$
at the end)
$$
\Big(f(x_{00}),f^{]1[}(x_{00},x_{01},t_{01}),f^{]1[}(x_{00},x_{10},t_{10}),
f^{]2[}((x_{00},x_{01},t_{01}),(x_{10},x_{11},t_{11}),t_{10}), \ttt \Big).
$$
For general $k$, as in \cite{Be05}, we will use multi-indices 
$\alpha \in  I_k := \{ 0,1 \}^k$,
namely $f^{[k]}$ depends on $2^k$
``space variables'' $x_\alpha$ with $\alpha \in I_k$ and on $2^k - 1$ 
``multitime variables''
$t_\alpha$ with $\alpha \in I_k^*:=I_k \setminus \{ 0 \}$.
We may thus denote the values of the difference quotient functions by
$$
f^{[k]}(\xxx,\ttt), \quad \quad \Delta^{[k]}(\xxx,\ttt), \quad \quad
\xxx \in V^{2^k}, \, \ttt \in \K^{2^k - 1}.
\eqno (1.1)
$$
It is then not too difficult to find the formula expressing $\Delta^{[k]}f$
in terms of all $f^{[j]}$, $j=0,1,\ldots,k$, namely, for $\alpha \in I_k$,
the $\alpha$-component of $\Delta^{[k]} f(\xxx,\ttt)$ is simply $t_\alpha$
for the ``multi-time variable'', and for the
$\alpha$-component of the ``space variable'' we get
$$
f^{[ \ell ]} \big( (x_\beta)_{\beta \subseteq \alpha},
(t_\beta)_{\beta \subseteq \alpha \atop \beta \not= 0} \big)
$$
where $\ell = |\alpha| = \sum_i \alpha_i$ is the ``depth'' of $\alpha$
and $\beta \subseteq \alpha$ means that $\beta_i \leq \alpha_i$ for all $i$
(this gives us $2^\ell$ choices for $\beta$, which is exactly the number
of arguments needed). In particular, the component with maximal depth $k$
equals $f^{[k]}(\xxx,\ttt)$. 
We leave the proof as an exercise to the reader.
However, the ``fine tuning'' of notation and formulae (our notation so far
does not impose a ``canonical order'' for the arguments!) is likely to depend on the
solution of the problems related to {\em scalar extensions}, see
Problem 5 below.

\msk
Back to first order difference calculus ($k=1$):
the reader may have already remarked that, as long as an invertible scalar
$t$ is fixed, difference calculus with respect to $t$ is rather trivial:
via the (linear) change of variables $(y,w):=(x,x+tv)$, i.e., via conjugating with
the matrix 
$\begin{pmatrix} 1 & 0 \cr 1 & t \cr \end{pmatrix}$,
the difference quotient map
$(x,v) \mapsto (f(x),\frac{f(x+tv)-f(x)}{t})$
is simply conjugate to $f \times f$,
and iterating this procedure for invertible values of $t$ essentially gives 
iterated direct products of $f$ with itself.
However, this change of variables becomes singular as soon as $t$ becomes
non-invertible, and hence this trivial interpretation of difference calculus
does {\em not} extend to differential calculus.
In other words, this kind of ``structure'' of difference calculus does not
(or at least not directly) solve our Problem 2, and it is important
to keep track of the dependence on $t$ and to pay attention whether this
dependence extends to singular values. In order make clear that we want
to include such singular parameter values in our considerations, we
will continue our list of problems
below  under the headline
``differential problems'', and we will use systematically the superscripts
$[k]$ instead of $]k[$.

\msk
Another comment may be helpful here:
in the case of functions of {\em one} variable, i.e.\ of {\em
curves} $f:I \to W$,  where
$I \subset \K$ is open, one may define inductively another kind of
difference quotients, called {\em divided differences}, by putting
$$
f^{>1<}(r,s)=\frac{f(r) - f(s)}{r - s }
$$
and
$$
f^{>k+1<}(t_1,\ldots,t_{k+1}):= {f^{>k<}(t_1,\ldots,t_k) -
f^{>k<}(t_2,\ldots,t_{k+1}) \over t_1 - t_{k+1}}.
$$
Then one easily sees that this expression is symmetric in its $k+1$ arguments,
and one can prove the explicit formula (see Chapter 7 in \cite{BGN04},
or \cite{Sch},  \cite{Rob} or \cite{Enc}):
$$
f^{>k<}(t_1,\ldots,t_{k+1}) = \sum_{j=1}^{k+1}
{ f(t_j) \over \prod_{i \not= j} (t_j - t_i) }
$$
If $f$ is $\CC^k$ in the sense defined above, then the function
$f^{>k<}$ admits a continuous extension to a continuous map $f^{<k>}$
defined on $I^{k+1}$  (including   the ``singular set'', where some of the
$t_i$'s coincide), and the $k$-th derivative $f^{(k)}(t)$ satisfies
$$
k! f^{(k)}(t)= f^{<k>}(t,\ldots,t). \eqno (1.2)
$$
In fact, the map $f^{<k>}$ can be seen as the restriction of $f^{[k]}$
to a suitable subset of $I^{[k]}$.
The converse of this statement is considerably more difficult:
if $\K$ a field and not merely a ring, and
if $f^{>k<}$ admits a continuous extension onto $I^{k+1}$, then
$f$ is $\CC^k$ (see \cite{BGN04} or \cite{Be08}, Exercices B.6).

\bigskip \nin
{\bf Problem 3: $q$-difference calculus.}
{\em 
Already on a formal level, it would be much more  satisfying if the number
of arguments for  $f^{[k]}$ and for
$\Delta^{[k]}f$ were a power of $2$, and not the odd number
$2^{k+1} - 1$.
In other words, it seems as if one ``time-like'' parameter were missing.
But what should be its meaning?
We ask for a ``correct'' definition of
the difference quotient maps, depending on an additional variable $t_0$ or $q
\in \K$.

\ssk More specifically, we have the impression that this additional variable
 somehow should play  a similar r\^ole as the ``quantum'' parameter
$q$ in the notion of {\em $q$-derivatives} (cf.\  {\em \cite{KC}}).
In particular, for $q=1$ one should find  difference and differential calculus
as introduced above, and for $q=-1$ one would like to find some sort of
super-difference and super-differential calculus, see Problem 9  below.}

\msk \nin {\bf Comments.} 
Here is a tentative definition, which is likely to be not yet the right guess:
recall from above the ``depth'' $\ell = |\alpha|$ of a variable 
$x_\alpha$, resp.\ $t_\alpha$. Starting from the definition of 
 $f^{[k]}$ and  $\Delta^{[k]} f$ as above, we may define functions
$\tilde f^{[k]}$, resp.\ $\tilde \Delta^{[k]} f$ by adding $t_0$ as new
variable (of depth $0$), and by replacing each argument $x_\alpha$
by $t_0^{|\alpha|} x_\alpha$, resp.\  $t_\alpha$
by $t_0^{|\alpha|} t_\alpha$, and finally re-dividing the resulting
expression again by $t_0^j$, where $j$ is the depth of the component in question.
For $t_0 = 1$, we get the old expressions.
For $t_0 = -1$, this definition means that we
 let act $\K^\times$ as if $V$ were replaced by its dual
space $V^*$: the natural action of the dilations on $V^*$ is the ``contragredient
representation''
$$
\K^\times \times V^* \to V^*, \quad (t,\phi) \mapsto \phi \circ (t \, \id)^{-1} =
t^{-1} \phi.
$$
My impression is that, on a foundational level, both actions of
$\K^\times$ on $V$ should have ``equal rights'', and the parameter $t_0$
should take care of this. It is obvious that the parity of
depth induces a $\Z/2 \Z$-grading on the space of our variables, and hence
we are automatically in a ``super-space'' setting -- compare with
 Problem 9 below.
%
%

\bigskip \nin
{\bf Problem 4: Barycentric difference calculus.}
{\em Instead of increasing the already incredibly big number of variables,
it might seem more reasonable to try to reduce this number as much as possible.
Indeed, one has the impression that much of the information encoded in the
maps $f^{[k]}$ is redundant, and that it should be possible to concentrate
this information in a map having less variables.
We have remarked above that for {\em curves} (maps of one scalar variable)
this is indeed possible, by looking at the ``divided differences''
$f^{<k>}$, which are maps of $k+1$ variables.
Is it possible to do something similar in the general case?}

\msk \nin {\bf Comments.}
My impression is ``yes, this cane be done''.
It should be achieved by some more subtle version of the trivial change of variables
$(x,v) \leftrightarrow (y,w)=(x,x + tv)$ mentioned above -- at first order, this change
becomes singular for non-invertible $t$, but at higher order only the
invertibility of $t_\alpha$ of depth $|\alpha|=1$ really matters; it seems that
 there ought to be a good change of variables involving all $t_\alpha$ with
depth $|\alpha| > 1$,  remaining continuous when these $t_\alpha$ tend to zero,
and such that the structure after change of variables becomes similar to the
one for the divided differences of a curve.

\ssk
Geometrically, this change of variables should correspond to changing from
the ``vector space point of view'' (the special r\^ole of $x_0$ corresponding to
the r\^ole of the origin in a tangent space) to the ``affine space point of view'':
we have to rewrite difference quotients by using barycentric calculus.
For instance, one could define inductively some kind of divided differences by
$$
f^{>1<}(x,y;s,t):= {f\big(sx+(1-s)y\big) - f\big(tx + (1-t)y \big) 
\over s-t } \quad \quad
(x,y \in V, s,t \in \K),
$$
whenever this makes sense; this map clearly contains the same information as the
old $f^{]1[}$. Then one would like to define inductively $f^{>k<}$ by a
recursion formula as a map 
of $2k+2$ variables (half of them space and half of them time variables),
similar to the recursion formula for the case of curves.
The ``explicit formula'' should be some weighted sum of values of $f$,
evaluated at barycenters defined in terms of the $k+1$ space variables
and the $k+1$ time variables.
As in the case of curves, the higher order differentials will then appear
in their ``divided form'' (that is, if we multiply by the factor $k!$, we get
the usual differentials, see the formula for curves above).
Geometrically, this corresponds to the ratio of the volume of
cubes (usual vector calculus) and simplices (barycentric calculus).
Both points of view (vectorial and affine) are important; in usual multivariable
calculus the vectorial point of view dominates, but it might turn out that many
of the more subtle problems could be better discussed via the affine (barycentric)
approach (see, for instance,  Problem 7 below).
Summing up, it seems as if in this issue
something very fundamental were going on.

\section{Differential problems}

As explained above, it is important to understand the dependence 
of difference quotients on the  ``multi-time'' parameter
$\ttt$.
On the other hand, we also have to study the behavior of these expression for
fixed  $\ttt$ (in such a way that the results remain valid for
singular $\ttt$). It seems to me that, for fixed $\ttt$, the point of view of
{\em scalar extensions} is most suitable, see the following

\msk \nin {\bf Exercise.}
{\em Fix a value $\ttt := (t_\alpha)_{\alpha \in I_k^*}$ for the ``multi-time parameter''
and let us consider the partial maps of $\Delta^{[k]} f$ for fixed $\ttt$, i.e.,
define
$$
\Delta_\ttt^{[k]} U  :=  \{ \xxx \in V^{2^k} | \, (\xxx,\ttt) \in \Delta^{[k]} U \},
$$
$$
\Delta_\ttt^{[k]} f  :  \Delta_\ttt^{[k]} U \to W^{2^k},  \quad  \xxx  \mapsto 
\pr_1 (\Delta^{[k]}  (\xxx,\ttt)) ,
$$
where $\pr_1(\xxx,\ttt)=\xxx$ is the projection onto the spacial part of the variable
$(\xxx,\ttt)$.
Show that $\Delta_\ttt^{[k]}$ is a covariant functor commuting with direct
products, and 
 that this functor can be interpreted as  the functor of scalar extension
from the ring $\K$ to a (commutative unital) ring $\Delta_\ttt^{[k]} \K$.
(In particular, if $f$ is a polynomial map, then
$\Delta^{[k]}_\ttt f$ is its scalar extension in the algebraic sense.)
}

\msk \nin {\bf Hints.}
Let us start with $k=1$, so we look at
$$
\Delta_t^{[1]} f (x,v)=\bigl(f(x),f^{[1]}(x,v,t) \bigr) = \bigl( f(x),
{f(x+tv)-f(v) \over t} \bigr),
$$
the latter provided $t \in \K^\times$. For such $t$,
the functoriality is easily checked (see Equation (0.4) above), and so is the
property $\Delta_t^{]1[} (f \times g)=\Delta_t^{]1[} f \times \Delta_t^{]1[} g$
(under obvious identifications of sets).
%
Again by density, the corresponding statements remain true for non-invertible 
scalars $t$; in particular, for $t=0$ we get again the chain rule for the
 {\em tangent functor}
$T = \Delta_0^{[1]}$, as defined by Equation (0.6).

\ssk
Now, for any $t \in \K$, the functor $\Delta_t^{[1]}$, being  covariant and
commuting with direct products, may be applied to  the ring $\K$ and its structure
maps $\K \times \K \to \K$,  and then
gives a new ring, which has dimension 2 over $\K$. For invertible $t$, 
a straightforward calculation gives the ring multiplication
\begin{eqnarray*}
(x_0,x_1)  \cdot (y_0,y_1)& =& 
\bigl(x_0y_0 \, , \,
{(x_0 + t x_1)(y_0 + t y_1) - x_0 y_0 \over t} \bigr) \cr
&  = & \bigl( x_0y_0 \, , \,  x_0y_1 + x_1y_0 + t x_1 y_1 \bigr) . \cr
\end{eqnarray*}
In a similar way, we see that the sum in this ring
 is just the usual sum in $\K^2$. Hence
as a ring, we get $\K \oplus \omega \K$ with relation $\omega^2 = t \omega$.
It can also be described as the truncated polynomial ring $\K[X]/(X^2 - t X)$.
Again by density, these statements remain true for non-invertible scalars $t$,
and in particular for $t=0$ we obtain the  {\em tangent ring} $T\K$,
which is nothing but the ring of  {\em dual numbers over $\K$}, $\K[X]/(X^2)=
\K \oplus \eps \K$
with $\eps^2 = 0$.

\ssk
Next, we claim that, if $f:U \to W$ is $\CC^\infty$ over $\K$,
then $\Delta_t^{[1]} f : \Delta_t^{[1]} U \to \Delta_t^{[1]} W$ is 
$\CC^\infty$ over the ring $\Delta_t^{[1]} \K$.
In fact, here the proof of the special case $t=0$ from \cite{Be05}, Theorem 6.2,
can be applied word by word; it uses only the fact that
$\Delta^{[1]}_t$ is a covariant functor commuting with direct products and
with diagonal mappings.  
What comes out is this: 
$$
\Delta^{[2]}_{t_{01},t_{10},t_{11}} f =
\Delta^{[1]}_{t_{01} + \omega t_{11}} (\Delta_{t_{10}}^{[1]} f) 
\eqno (2.1)
$$
where $\omega$ satisfies $\omega^2 = t_{10} \omega$.
Thus $\Delta^{[2]}_\ttt$ is a composition of two scalar extension functors.
 By induction, it  follows now that  $\Delta_\ttt^{[k]}$ is a composition
of $k$ scalar extension functors, and hence is itself a 
functor of scalar extension from $\K$ to a ring $\Delta_\ttt^{[k]}\K$.

\msk \nin
{\bf Problem 5: Scalar extensions.}
{\em Describe  the structure of the  ring $\Delta_\ttt^{[k]}\K$ in such a
way that its dependence on $\ttt$ becomes understandable.
}

\msk \nin {\bf Comments.}
For $k=1$, as seen above, $\Delta^{[1]}_t\K \cong \K[X]/(X^2 - t X)$.
For $k=2$, by (2.1), this ring is isomorphic to
$$
\Big(\K[X_1]/(X_1^2 - t_{10} X_1)\Big)[X_2]
/(X_2^2 - (t_{01} + t_{11} X_1) X_2) 
$$
$$
\cong
\K[X_1,X_2]/\big( (X_1^2-t_{10} X_1), (X_2^2 - (t_{01} + t_{11} X_1) X_2)
\big).
$$
We can also describe this ring as
$$
R = \K \oplus \omega_{10} \K \oplus \omega_{01} \K  \oplus \omega_{11} \K
$$
with relations
$$
\omega_{10}^2 = t_{10} \omega_{10}, \quad
\omega_{11} = \omega_{01} \omega_{10} = \omega_{10} \omega_{01}, \quad
\omega_{01}^2 = (t_{01} + t_{11} \omega_{10}) \omega_{01} =
t_{01} \omega_{01} + t_{11} \omega_{11}.
$$
It follows that
$$
\omega_{10}\omega_{11}=t_{10}\omega_{11}, \quad
\omega_{01}\omega_{11}=(t_{01}+t_{10}t_{11})\omega_{11},
$$
$$
\omega_{11}^2 = \omega_{10}^2 \omega_{01}^2 = 
 t_{10} \omega_{10}(t_{01} + t_{11} \omega_{10}) \omega_{01} =
 t_{01} (t_{10}+t_{01}t_{11}) \omega_{11}.
$$
Some features of
the structure of this ring are similar to the ones
 of the second order tangent ring $TT\K$
(cf.\ \cite{Be05}).
For general $k$, we may write
$$
R = \bigoplus_{\alpha \in I_k} \omega_\alpha \K \quad \quad (\omega_0 = 1),
$$
with relations
$$
\omega_\alpha \omega_\beta = 
\sum_\gamma c_{\alpha,\beta}^\gamma \omega_\gamma ,
$$
where the coefficients
 $c_{\alpha,\beta}^\gamma$ depend on the multi-time parameter $\ttt$.
It seems that this system of coefficients is ``triangular'' in a certain
sense, and there should be some general pattern permitting to understand it.

\bigskip \nin
{\bf Problem 6: Difference version of the de Rham complex.}
{\em We ask for a (first) necessary condition for solving difference 
equations: 
 is there an operator $D$ whose square is zero and such
that $DF=0$ is a necessary condition for $F:U^{]k[}\to W$ to be 
essentially of the form $f^{]1[}$ for some $f:U^{]k-1[} \to W$?
}

\msk \nin
{\bf Comments.}
The classical de Rham complex concerns the case $\ttt = 0$, and its
construction relies entirely on the fact that the symmetric group
$\Sigma_k$ acts by automorphisms of the ring $T^k \K$ (or, equivalently,
of the functor $\Delta_{\bf 0}^{[k]}$), cf.\ \cite{Be05}, Chapter 22.
Therefore also for general $\ttt$ one has to study the ``Galois group''
$\Aut_\K(\Delta^k_\ttt \K)$.
It should be some kind of deformation of the automorphism group
$\Aut_\K(T^k \K)$, and we expect that its quotient with respect to
its ``identity component'' is still isomorphic to $\Sigma_k$.
Therefore we should indeed have necessary conditions in terms of this
group for a map $U^{[k]} \to W$ to be some $f^{[1]}$, and these conditions
 should look similar to closedness of $k$-forms.

\bigskip \nin
{\bf Problem 7: The case of positive characteristic.}
{\it 
Differential calculus and differential geometry work perfectly well over
base fields of arbitrary, possibly {\em positive}, characteristic
(see {\em \cite{BGN04}, \cite{Be05}}). 
However, in case of characteristic $p>0$,
it is less obvious how the ``differential information'' is best
exploited. The approach from   {\em \cite{Be05}}, based on the higher order
tangent functor $T^k$ and on its invariance under the symmetric group
(see the preceding problem), still leads to some loss of information.
How can we avoid this loss of information?}

\msk \nin
{\bf Comments.}  It is known that, in case of characteristic $p>0$,
 the usual definition of a {\em Lie algebra} is too weak, and that one
should use stronger structures such as the {\em restricted Lie algebras}
or {\em Lie $p$-algebras} introduced by N.\ Jacobson \cite{Jac41} (and, later,
for related algebras, such as {\em Leibniz algebras}, cf.\ \cite{DL}).
Roughly, one has to add to the usual Lie algebra structure a compatible structure
of a ``$p$-th power map'' $x \mapsto x^{[p]}$. 

\ssk
It seems as if for differential geometry in general (which comprises Lie
theory) the situation were very similar: when the approach is based
on usual higher differentials $d^k f$ (or, equivalently, on the higher order
tangent functors $T^k =\Delta^{[k]}_{\ttt}$ for $\ttt =0$), 
we inevitably lose information;
but this loss of information can be avoided if we use the information
encoded in the higher order slopes $\Delta^{[k]}_\ttt$ 
where {\em some, but not all} $t_\alpha$ are zero. 
In fact, the general Taylor formula, valid in any characteristic,
expresses the Taylor coefficients $a_j$ in the expansion
$$
f(x+th)=f(x)+\sum_{j=1}^k t^j a_j(x,h) + t^k R(x,h,t)
$$
in terms of the slopes (see \cite{BGN04}, \cite{Be05} or \cite{Be08}):
$$
a_1(x,h)=f^{[1]}(x,h,0), \quad a_2(x,h)=f^{[2]}((x,h,0),(0,0,1),0) , \ldots
$$
In case of a curve $f:I \to W$, this can also be expressed in the context 
of divided differences (cf.\ Problem 2):
$$
a_1(x,h)=f^{<1>}(x) \, h, \quad a_2(x,h)=f^{<2>}(x,x)\, h^2 , \ldots ,
a_k(x,h)=f^{<k>}(x,\ldots,x)\, h^k.
$$
For instance, in characteristic 2, the second differential of $f(x)=x^2$ 
at the origin is zero, but nevertheless we  get the correct Taylor expansion:
we have $f^{<1>}(t,s)=\frac{s^2-t^2}{s-t} = s+t$ and
$f^{<2>}(r,t,s)=\frac{(r+t)-(r+s)}{t-s}=1$, whence the Taylor expansion
$h^2 = f^{<2>}(0,0,0)h^2$. In fact, in the general case, it is easily
checked that
 $f^{[2]}((x,h,t),(0,0,1),s))$  is exactly the 
divided difference $\gamma^{<2>}(0,t,s)$ 
for the curve $\gamma(t)=f(x+th)$, thus
relating the Taylor expansion for the curve $\gamma$ at $0$
 and the map $f$ at $x$.
Put differently, the (hypothetic) ``barycentric'' or ``simplicial
 differential calculus'' from
Problem 4 should furnish the information that is lost in characteristic $p>0$
if we use the classical ``vectorial'' or ``rectangular differential calculus''.

\ssk
It is all the more surprising that, even in positive characteristic, one can
{\it reconstruct} some sort of approximation of a Lie group by starting from
mere Lie- or Leibniz algebras (see \cite{Di07}); it looks as if one just had
to add some structure like the $p$-th order power map in order to fill in
the missing information.

\bigskip \nin
{\bf Problem 8: Non-commutative base rings.}
{\it 
In principle, the definition of the maps $f^{[k]}_\K$ and
$\Delta^{[k]}_\K f$
 makes sense over non-commutative base fields or rings $\K$.
However, already the multiplication map $\K \times \K \to \K$ and the 
squaring map $\K \to \K$ will then no longer be of class $\CC^1$, and hence
all arguments relying on these facts break down.
How can definitions be modified in order to take account of these problems?}

\msk \nin
{\bf Comments.} If $m:\K \times \K \to \K$ is the product in the ring $\K$, then
$$
m^{]1[}((x,v),(y,w),t)= t^{-1} \big( (x+tv)(y+tw) - xy \big) =
vtw +  vy + t^{-1} xtw .
$$
Unfortunately, the last term does in general not admit a  continuous extension
to $t=0$ if $\K$ is not commutative. Thus $m$ will not be of class $\CC^1$
over $\K$,
and similarly, the power maps $\K \to \K$, $x \mapsto x^k$ won't be $\CC^1$ over
$\K$
either. Therefore we do not get a useful differential calculus over non-commutative
fields or rings.
(Of course, one could define the difference quotient map by putting
$t^{-1}$ to the right of the bracket, but this does not help!)

\ssk
There is no simple way to avoid this problem -- as is well-known,
 no quaternionic calculus has ever been found that in any respect
shared the power of real or complex calculus. However, it is not excluded
to generalize some more specific aspects of calculus to the non-commutative case.
To do this, one certainly needs to impose some suitable conditions on $\K$.
For instance, one might rather consider {\em Hermitian rings}, that is
a pair $(\K,\tau)$ formed by a ring together with an involution (anti-automorphism
of order 2). One may then try to use the involution  for an appropriate definition
of difference quotients. This might be interesting in relation with exceptional
Lie groups and their geometries. (It is known that all compact
forms of exceptional Lie groups arise as automorphism groups of {\em quaternionic
symmetric spaces}, the so-called Wolf spaces. These spaces  are not really manifolds
over the quaternions, but it looks as if they belonged to some sort of generalized
projective geometries (in the sense of \cite{Be02})
 that should be defined in some suitable way
 over the quaternions.)

\bigskip \nin
{\bf Problem 9: Super differential calculus.}
{\em Is it possible to see {\em super differential calculus}  as a continuous extension of some ``super difference calculus''
to singular values?}

\msk \nin {\bf Comments.}
In the context of Problem 8,
another suitable condition to impose on our ring $\K$ would be 
{\em super-commutativity}, i.e., $\K = \K_1 \oplus \K_{-1}$ is
$\Z /2\Z$-graded and the product satisfies the sign rule
$xy=(-1)^{d(x) \cdot d(y)}yx$.  If it is possible to define some
reasonable difference calculus in this case, one should expect it to have
some link with {\em super differential calculus} 
(cf.\  \cite{Dewitt}, \cite{Var}).

\ssk
Another aspect of this problem is given by
 the scalar extension point of view: in presence of a
``quantum'' parameter $q$ (see Problem 3), the arguments used above (hints to
the exercise) 
suggest that for $q=-1$ we should obtain some
calculus that can be interpreted as a scalar extension by the Grassmann algebra.

\bigskip \nin
{\bf Problem 10: Pansu calculus.}
{\em Another approach to non-commutativity comes from sub-Riemannian geometry:
instead of making $\K$ non-commutative we can ``make the vector space
 $V$ non-commutative'' (the co-called {\em Pansu calculus})\footnote{I thank
Marius Buliga for pointing out to me this approach, cf.\ \cite{Buliga}.}. 
Is there a corresponding ``Pansu difference calculus'', having the
Pansu differential calculus as limit, and if so, what are its
properties?}

\msk \nin {\bf Comments.}
The non-commutative version of a vector space is a {\em Carnot group}: let
$V = \g = \g_1 \oplus \ldots \oplus \g_n$ be a graded nilpotent Lie algebra
(i.e., $[\g_i,\g_j] \subset \g_{i+j}$; wo do not require here equality as in
\cite{Buliga}) and
define a (polynomial) group law $v * w = v+w+\frac{1}{2}[v,w]+\ldots$ 
on $V$ by the Campbell-Hausdorff formula.
If $\g$ is abelian, then this is just usual vector addition, and hence
$(V,*)$ can be seen as a non commutative analog of $(V,+)$.
Now we define the analog of multiplication by scalars:
for $t \in \K$ define the {\em dilation}
$$
\delta_t : \g \to \g, \quad \sum_i x_i \mapsto \sum_i  t^i x_i.
$$
From the grading condition it follows immediately that $\delta_t$ is a
Lie algebra-endomorphism
of $\g$; if $t$ is invertible, it is an automorphism.
If $\g=\g_1$, we get the usual multiplication by scalars.
Now let $f:U \to W$ be a map between Carnot groups and
define the {\em first difference quotient} by
$$
f^{]1[}(x,v,t):= \delta_t^{-1} \Big( f(x)^{-1} * f(x * (\delta_t v) ) \Big).
$$
(Of course, since $(V,*)$ is no longer commutative, we had to make a choice
of order:
all other choices have equal rights!)
Now assume that our Carnot groups are topological groups, that
the dilation map $\K \times V \to V$ is continuous and that $U$ is open
in $V$. Then we can define the class $\CC^1$ as before by requiring that
$f^{]1[}$ extends continuously to the suitably defined set $U^{[1]}$.
To some extent, the theory goes through: for instance, one proves as before
that the ``differential'' $df(x):V \to W$ is a ``Carnot linear map''.
Things start to become more involved when looking at higher order different quotients
and trying to define the class $\CC^2$:
before iterating our procedure, we have to fix a Carnot structure on direct
products; the direct product being $\Z \times \Z$-graded, there are several
possible choices of $\Z$-gradings. It seems reasonable to choose the
``diagonal grading''; but even then we do not get a well-defined notion of
``bilinear maps'', and hence we have to be very careful when talking about
higher order differentials.
Nevertheless, this sort of calculus might be interesting for geometries
such as certain flag varieties that are indeed ``Carnot groups patched together''
and which could be smooth manifolds in the sense of such a calculus.

\bigskip \nin
{\bf Problem 11: Pointwise concepts.}
{\em Can one define the notion of a map $f$ being ``of class $\CC^k$ at
one point''? Is it then true that $f$ is of class $\CC^k$ on $U$ if $f$ is $\CC^k$ at all
points $a \in U$?}

\msk \nin {\bf Comments.}
This should indeed be possible, but as far as I know it has not yet been
really worked out. For $k=1$ and $k=2$, I have outlined the definition of the
``classes $\CC^k_a$'' in \cite{Be08}, Section 14.3:
for $k=1$, we require that $f$ is continuous at $a$ and that
the difference quotient $f^{]1[}$ admits an extension, denoted by
$f_a^{[1]}$, to 
$$
U_a^{[1]}:= U^{]1[} \cup V_a, \quad \mbox{where} \quad
V_a := \{ (a,v,0) | \, v \in V \},
$$
such that $f_a^{[1]}$ is continous on $V_a$.
This notion has good properties, and $f$ is indeed $\CC^1$ on $U$ iff it is
$\CC^1_a$ for all $a \in U$. 
In finite dimension over $\K=\R$, this property is equivalent to the
{\em strict differentiability of $f$ at $a$} (see loc cit). It is known that strict
differentiability is in many respects the ``good pointwise concept'' (better
than the weaker Fr\'echet differentiability); in particular, $f$ is
$\CC^1$ on $U$ iff it is strictly differentiable at all points of $U$.
The definition of the class $\CC^k_a$ is similar, though slightly more subtle:
roughly, we require the existence of an extension of $f^{]k[}$ onto 
a set $U_a^{[k]}$ defined similarly as above, and requiring continuity on
the set of all multi-space variables $x_\alpha$ of depth $|\alpha| \geq 1$ and
all multi-time variables $t_\alpha$ of depth $|\alpha| \geq 2$.
This still gives a ``pointwise'' (as opposed to ``local'') concept having
properties strong enough to prove all basic results of differential calculus.

\bigskip \nin
{\bf Problem 12: Integral formulae in the real case.}
{\em In the case $\K=\R$, when dimensions are finite or in presence of a
locally convex topological vector space structure, the fundamental theorem
of calculus permits to express first order
difference quotients by integrals 
(cf.\  \cite{BGN04} for the general locally convex case):
$$
f^{[1]}(x,v,t)=\int_0^1 df(x+stv)v \, ds.
$$
Find an analoguous  expression of $f^{[k]}(\xxx,\ttt)$ by an integral formula !}

\msk \nin {\bf Comments.}
In principle, we just have to iterate the preceding formula:
$$
f^{[2]}((x,v,t),(x',v',t'),s) = \int_0^1 d(f^{[1]})((x,v,t)+s' 
(x',v',s')) (x',v',t') \, ds' \, ,
$$
then again express $f^{[1]}$ by an integral and permute the differential
$d$ with integration. We will end up with a certain integral of 
$d^2 f$ over a plane region, and so on. 
Via Problem 4, this question is related to the corresponding question
for divided differences of
curves. In this case, on the following integral formula holds
$$
f^{<k>}(x_0,\ldots,x_k)= \int_{\Delta_k} f^{(k)}\bigl( x_0 +\sum_{i=1}^k
\lambda_i (x_i-x_0) \bigr) \, d\lambda_1 \cdots d\lambda_k \, 
$$
(see, for instance, \cite{Be08}, Exericise B.6.7), where 
$f^{(k)}$ is the usual $k$-th order derivative of a curve and
$$
\Delta_k = \{ (\lambda_1,\ldots,\lambda_k) \in \R^k \vert \,
\lambda_i \geq 0, \sum_i \lambda_i \leq 1 \}
$$ 
is the standard simplex in $\R^k$. (Recall that its volume is $\frac{1}{k!}$).
By averaging further, one deduces that there exists a function
 $M_{k+1}(x_0,\ldots,x_k;t)$ (called the {\it Peano kernel}), such
that, for all curves of class  $\CC^{k}$, $f:I \to \R^n$,
$$
f^{<k>}(x_0,\ldots,x_{k}) = \frac{1}{k!} \int_\R
f^{(k)} (t) \cdot M_{k+1}(x_0,\ldots,x_k;t) \, dt .
$$
One can show that $M_{k+1}(x_0,\ldots,x_k;\cdot)$ is a function of
class $\CC^{k-2}$ such that its  restriction to the interval
$]x_j,x_{j+1}[$ is a polynomial of degree $\leq k-2$, and zero
outside the interval $]x_0,x_k[$.
%

\bigskip \nin
{\bf Final remarks.}
In some sense, all the problems mentioned so far are only preliminaries
for the Big Problem which is the {\em integration problem}:
{\em is there a difference theory of an anti-derivative?}
I am inclined to believe that such a theory exists, but it is likely to be
fairly sophisticated, rather giving an infinite sequence of necessary
conditions than a plain ``existence and uniqueness theorem''.
In any case, before understanding anti-derivatives
and ``anti-differences'', we have to understand derivatives and differences.


}

\vfill \eject

\end{document}